\documentclass{amsart}
\usepackage{amssymb,stmaryrd}
\usepackage{amsfonts}
\usepackage{amstext}
\usepackage{algorithmic}
\usepackage{algorithm}
\usepackage{graphicx}
\usepackage{epstopdf}
\usepackage[all]{xy}

\numberwithin{equation}{section}
\newtheorem{theorem}{Theorem}[section]
\newtheorem{lemma}[theorem]{Lemma}
\newtheorem{proposition}[theorem]{Proposition}

\theoremstyle{definition}
\newtheorem{definition}[theorem]{Definition}
\newtheorem{example}[theorem]{Example}

\theoremstyle{remark}

\newcommand{\CA}{{\mathcal{A}}}
\newcommand{\C}{{\mathbb{C}}}
\newcommand{\Z}{{\mathbb{Z}}}

\newcommand{\N}{{\mathbb{N}}}

\newcommand{\CC}{{\mathcal{C}}}
\newcommand{\CH}{{\mathcal{H}}}

\newcommand{\CR}{{\mathcal{R}}}
\newcommand{\CS}{{\mathcal{S}}}

\newcommand{\CM}{{\mathcal{M}}}
\newcommand{\CZ}{{\mathcal{Z}}}

\renewcommand{\ker}{{\rm{ker}}}

\newcommand{\isom}{\cong}

\newcommand{\tens}{\otimes}

\newcommand{\id}{{\rm id}}
\newcommand{\bo}{{}^{(1)}}
\newcommand{\bt}{{}^{(2)}}
\renewcommand{\o}{{}_{(1)}}
\renewcommand{\t}{{}_{(2)}}
\renewcommand{\th}{{}_{(3)}}

\newcommand{\eps}{\epsilon}

\newcommand{\und}{\underline}

\newcommand{\la}{{\triangleright}}
\newcommand{\ra}{{\triangleleft}}
\newcommand{\lbiprod}{{>\!\!\!\triangleleft\kern-.33em\cdot}}
\newcommand{\rbiprod}{{\cdot\kern-.33em\triangleright\!\!\!<}}
\newcommand{\lcross}{{>\!\!\!\triangleleft}}

\newcommand{\lrbicross}{{\blacktriangleright\!\!\!\triangleleft}}

\newcommand{\codcross}{{\blacktriangleright\!\!\blacktriangleleft}}

\begin{document}

\title{Strict quantum 2-groups}


\author{S. Majid}
\address{Queen Mary, University of London\\
School of Mathematics, Mile End Rd, London E1 4NS, UK}

\email{s.majid@qmul.ac.uk}


\begin{abstract}
A crossed module is $(A,H,d,\la)$ where $d:A\to H$ is a homomorphism of groups and $H$ acts on $A$, with conditions leading to a groupoid $A\lcross H{\to\atop \to}H$ as an example of a strict 2-group. We give the corresponding notion of a quantum 2-group where we replace the above by Hopf algebras and introduce a new version of quantum groupoid. The work also suggests a natural notion of braided crossed module where $A$  a braided-Hopf algebra in the braided category $\CZ({}_H\CM)$ of crossed $H$-modules, although without the full groupoid picture in this more general case.  \end{abstract}
\maketitle 

\section{Introduction} 

There has recently been much interest in higher algebra, motivated in part by the search for a better understanding of 4-dimensional topological invariants and quantum gravity in 4 dimensions as a generalisation of the successes in 3 dimensions. The latter has revolved around 3-manifold invariants constructed as Chern-Simmons theory via quantum groups\cite{Dri,Ma:book} as a view on 3-dimensional quantum gravity with point sources and generalisations have includes on the one hand 2- groups \cite{Bro,BL,BM} and 2-algebra and on the other the weak Hopf algebras, Hopf algebroids and related structures\cite{Tak,Sch,BrzM, Xu, Bo, BNS,BS}. In this paper we bring together two of these strands with a formulation of quantum 2-groups in line with that of 2-groups. This is very different both from existing 2-vector and 2-algebra ideas for 2-Hopf algebras (which are intimately tied up with the symmetric group). Rather, it is quantum-groupoid oriented but in order to make it work we need a new notion, which we call `embeddable quantum groupoid' (Definition~2.3) critically different from previous approaches in the case of a noncommutative base. This is then used in defining a 2-quantum group (Definition~2.5) and Theorem~2.6 provides a construction from natural Hopf algebra data. This is the main result of the paper. The crossed module data in our approach has a further generalisation using Hopf algebras in braided categories which appears very natural and this is given in Section~3. We only braid the data, which could have applications in their own right, without the braided generalisation of embeddable quantum groupoid which appears to require new ideas beyond those here.

We now explain the background of 2-groups \cite{Bro,BL,BM}. Just as a group $G$ may be regarded as a category with one object $*$ and morphisms $G$, a strict 2-group is an example of a 2-category with one object $*$,  morphisms $G_0$ and a set of 2-morphisms $G_1$. The latter are both groups used now for the morphism composition, all 1-morphisms  in $G_0$ have the same source and target $*$, while for 2-morphisms in $G_1$ there are source and target maps $s,t: G_1{\to\atop \to}G_0$ both of them group homomorphisms. We also require a compatible identity 2-morphism for every 1-morphism, i.e. an algebra map $i: G_0\to G_1$ with the obvious properties $s i=t i=\id$. In addition, $G_1$ has a second `vertical' composition $\circ$ of elements with the same source and target and an interchange law between vertical and horizontal compositions.

One can also think of such a strict 2-group as a strict monoidal category with objects $G_0$ and morphisms $G_1$ with the group products as the respective tensor products and $s,t$ the source and target maps for morphisms. There is also a second product $\circ$ on $G_1$ for the composition of morphisms and functoriality of the tensor product requires amounts to an  interchange law.  

A third, and related point of view is that a strict 2-group is a groupoid with $s,t$ the source and target maps and for which the two spaces are groups and $s,t,i$ are group homs. In addition the groupoid product $\circ$ obeys the interchange law $(wz)\circ(xy)=(w\circ x)(z\circ y)$ where $s(w)=t(x)$ and $s(z)=t(y)$ so as to be composable. In all three points of view we also require 1-morphisms to be `invertible' in a suitable sense where $x\circ x^{-}=i t(x)$, $x^{-}\circ x=i s(x)$ which we distinguish from the group inverse.

A simple example of such objects is provided by the following data: a group hom $d:A\to H$ between groups, an action $\la$ of $H$ on $A$ by group automorphisms, subject to 
\[ d(h\la a)=h d(a) h^{-1},\quad \forall a\in A,\ h\in H\]
(a precrossed module) and
\[ d(a)\la b= a b a^{-1},\quad \forall a,b\in A\]
(a crossed module). The action means that we have a cross product  
\[ A\hookrightarrow A\lcross H {{\buildrel s\over \to}\atop {\buildrel \hookleftarrow \over i}} H\]
with group homs given by the obvious coordinate-wise inclusion and projection $s(a,h)=h$, $i(h)=(e,h)$. $A$ appears as a normal subgroup here. We set $G_0=A$ and $G_1=A\lcross H$ while the data $d$ and the first axiom stated provide a map,
\[ t(a,h)=d(a)h\]
and the second more optional axiom evidently implies that the kernel of $d$ is a commutative subgroup. It is also reasonable for consistency with $d$ and the first condition. 

As a 2-category the group structures are the composition of 1-morphisms and horizontal composition of 2-morphisms respectively. The vertical composition of $(a,h): h\Rightarrow d(a)h$ followed by $(b,g):g\Rightarrow d(b)g$ is defined when $g=d(a)h$ and is given by $(ba,h):h\Rightarrow d(ba)h$ which makes sense as $d$ is a group hom. The interchange law then comes down to the last requirement above (giving a crossed rather than a precrossed module). Similarly, as a monoidal category the group structures are the $\tens$ and the composition of morphisms is
\[ (b, d(a) h)\circ (a,h)=(ba, h),\quad \forall a,b\in A,\ h\in H\]
where again the product is defined when the target of $(a,h)$ matches the source of the following morphism.  As a groupoid this $\circ$ is the groupoid structure. The inversion with respect to $\circ$ is
\[ (a,h)^{-}=(a^{-1},d(a)h),\quad \forall a\in A,\ h\in H.\]
Let us make a small observation here that  $\circ$ is actually defined associatively on all of $G_1\times G_1$ by $(b,g)\circ(a,h)=(ba,h)$ for all $a,b\in A$, $g,h\in H$ which restricts correctly on $G_1\square G_1\subseteq G_1\times G$ defined as the subspace where the source and target match up. This embedding of $\circ$  will be key to our approach.

We will make use of quantum group theory in the Hopf algebra approach of \cite{Ma:book}, and refer to this for details. Particularly, we use the Sweedler notation $\Delta h=h\o\tens h\t$ for the coproduct of the Hopf algebra (there is an implicit sum of terms here) and we require that our Hopf algebra $H$  has invertible antipode. We work over a field $k$. With care one may replace $k$ by $\N$ and perhaps eventually `categorify' these notions that way, but that is not something we consider here. Note that we use the term crossed module as above for the 2-group data and its generalistion. Closely related, we say $H$-crossed module as in \cite{Ma:book} for a Drinfeld-Radford-Yetter module over a Hopf algebra $H$. This is the same thing as an object of the centre $\CZ({}_H\CM)$ \cite{Ma:rep} of the category of $H$-modules, i.e. a module of the quantum double\cite{Dri} of $H$ appropriately defined. 

\subsubsection*{Acknowledgements} I would like to thank John Barrett for suggesting the problem. Also thanks to Behrang Noohi.

\section{Hopf crossed modules}

Our staring point is an immediate generalisation of the group crossed module data where we consider group algebras in place of groups and then let the former be arbitrary Hopf algebras:

\begin{definition} A Hopf algebra crossed module means:
\begin{enumerate}
\item a Hopf algebra $A$ in the category ${}_H\CM$ of modules of a Hopf algebra $H$ (i.e. an $H$-module algebra and coalgebra) such that
\[ h\o \tens h\t \la a = h\t \tens h\o\la a,\quad \forall a\in A,\ h\in H\]
\item a Hopf algebra map $d: A\to H$ such that
\[ d(h\la a)=h\o d(a) Sh\t,\quad \forall a\in A,\ h\in H\]
\item  
\[ d(a)\la b= a\o b S a\t,\quad \forall a,b\in A.\]
\end{enumerate}
\end{definition}

The conditions (1) and results in  \cite{Ma:phy} provide for the formation of a double cross product $\bowtie$, in our case with one action trivial, so we have a Hopf algebra $A\lcross H$. As an algebra it is the smash product
\[ (a\tens h)(b\tens g)=a (h\o\la b)\tens h\t g\]
and we have the tensor product unit and coalgebra. It is also a special case of a bicrossproduct $\lrbicross$ as in \cite{Ma:phy} with the coaction trivial. The smash product comes with the canonical Hopf algebra maps
\[ A\hookrightarrow A\lcross H {{\buildrel s\over \to}\atop {\buildrel \hookleftarrow \over i}} H,\quad i(h)=1\tens h,\ s(a\tens h)=\eps(a)h\]
given by inclusion via the relevant unit or projection by application of the relevant counit.  Next, the conditions (2) allow us to define the map
\[ t(a\tens h)=d(a)h,\quad\forall a\in A,\ h\in H,\]
which we will see later is necessarily a Hopf algebra map. The condition (3) then ensures that the subalgebra 
\[ A^{inv}=\{a\in A\ |\ d(a)=\eps(a)\}\]
is commutative.  Without this condition (3) we say that we have a Hopf precrossed module. Note also that (1) implies that $h\la Sa=S(h\la a)$ for all $h\in H$ and $a\in A$. Using this in (2) and (3) then gives significant constraints on $H$ in relation to the image of $d$ and on $A$, respectively.

We next give an interpretation of this construction which will make it a special case of the next section. We recall that if $H$ is any Hopf algebra the category $\CZ({}_H\CM)$ of crossed-$H$-modules (also called Drinfeld-Radford-Yetter modules) consists of vector spaces $V$ which are both (left) modules and comodules of $H$ such that
\[ \Delta_L (h\la v)= h\o v\bo S h\t \tens h\th\la v\bt\]
where $\Delta_L v=v\bo\tens v\bt$ is the coaction and $\la$ is the action. When $S$ is invertible this forms a braided category with braiding 
\[ \Psi_{V,W}(v\tens w)=v\bo\la w\tens v\bt,\quad \forall v\in V,\ w\in W.\]
Note also that  $H$ itself is an object in $\CZ({}_H\CM)$ under the adjoint action $h\la g=h\o g Sh\t$ and regular coaction $\Delta_L=\Delta$.

\begin{proposition} The conditions for a Hopf crossed module in  Definition~2.1 are equivalent to the following:

\begin{enumerate}
\item $A$ is a Hopf algebra in the category $\CZ({}_H\CM)$  with the given action and the trivial coaction $\Delta_La= 1\tens a$. 
\item There is a Hopf algebra map $d:A\to H$ such that $A$ becomes an object of $\CZ({}_H\CM)$ when viewed with the given action and the pushout coaction $\Delta_L=(d\tens\id)\Delta$ (then $A$ becomes an algebra in this category and $d$ a morphism)
\item The algebra $A$ is commutative in the category when viewed as in item (2). 
\end{enumerate}
\end{proposition}
\proof (1) That $A$ is an crossed $H$-module reduces in the case of  the trivial coaction to the condition displayed in (1) of Definition~2.1  for a given action $\la$. That $A$ is an algebra and coalgebra with respect to the $H$-module structure is that it is an $H$-module algebra and coalgebra as in (1) of Definition~2.1. That they do so with respect to the the trivial $H$-coaction imposes no conditions. Finally, that $A$ forms a Hopf algebra in the category with trivial coaction is immediate because the braiding of $\CZ({}_H\CM)$ on $A$ in this case is trivial, i.e. that assertion is equivalent to the assertion that it is an ordinary Hopf algebra.  
 
(2)  If we are given (1) and that $d$ a Hopf algebra map then $A$ becomes an $H$-comodule algebra by the pushout coaction $\Delta_L(a)=d(a\o)\tens a\t$. This then becomes an  crossed $H$-module with the given action,
\[ \Delta_L(h\la a):= d((h\la a)\o) \tens (h\la a)\t=d(h\o\la a\o)\tens h\t\la a\t=h\o d(a\o)Sh\t\tens h\th\la a\t\]
\[=h\o a\bo Sh\t\tens h\th\la a\bt=h\t a\bo Sh\th\tens h\o\la a\bt\]
as required provided the condition displayed in (2) of Definition~2.1 holds for the third equality. Conversely, if $d(h\o\la a\o)\tens h\t\la a\t=h\o d(a\o)Sh\t\tens h\th\la a\t$, which is the content of the crossed $H$-module compatibility for $\Delta_L$ given the other facts, apply $\eps$ to deduce the condition displayed in (2) of Definition~2.1. In this case $A$ becomes an algebra in $\CZ({}_H\CM)$ with this potentially new coaction $\Delta_L$ and we can check that $d$ is then a morphism: The part of this with respect to the action is the condition displayed in (2) of Definition~2.1 again, while  the part of this  with respect to the coaction is  $(\id\tens d)\Delta_L:=(d\tens d)\Delta=\Delta d$ which holds as $d$ is a coalgebra map. Our requirements in fact  force $d$ to be a Hopf algebra map if one wants to put the axioms in minimal form. 

(3) If the condition displayed in (3) of Definition~2.1 holds then 
\[ \Psi(a\tens b)=d(a\o)\la b\tens a\t=a\o b Sa\t \tens a\th\]
which in multiplying both sides says that $\cdot\Psi(a\tens b)=ab$. Conversely, if $(d(a\o)\la b)a\t=ab$ for all $a,b\in A$ then 
\[ d(a)\la b=(d(a\o)\la b)a\t Sa\th= a\o b S a\t\]
which is the condition displayed in (3) of Definition~2.1. So this is an equivalent formulation. \endproof

The conditions (1)-(2) define a Hopf precrossed module. We see that this is then a Hopf crossed module iff $A$ is commutative in the category $\CZ({}_H\CM)$ with respect to the induced coaction.

We now turn to exactly what kind of object can be built from this data. Let us note first that if $C$ is a coalgebra and $M,N$ are $C$-bicomodules (so equipped with mutually commuting left and right coactions) then 
\[ M\square N=\{\sum_i m_i\tens n_i\ |\ \sum_i \Delta_R m_i\tens n_i=\sum_i m_i\tens \Delta_Ln_i\}\subseteq M\tens N\]
(the `cotensor product') is another such, an operation which is associative up to the obvious identification. This is just the dual notion to the tensor products of bimodules over an algebra. This category ${}^C\CM^C$ of bicomodules with $\square$ is not, however, canonically braided when $C$ is not cocommutative, so we cannot define the notion of a cobialgebroid in the most obvious way as a bialgebra here. Although a notion of bialgebroid has been found in the literature and can be dualised, our example does not appear easily to fit within it and we shall formulate what seems more immediately to apply.
Thus for our purposes:

\begin{definition} An `embedded quantum groupoid' over a base coalgebra $C$ (over a field $k$) is:
\begin{enumerate}
\item a coalgebra $M$ over $k$, with coproduct $\Delta m=m\o\tens m\t$ and counit $\eps$, and coalgebra maps $s,t,i$,
\[ M {{{\buildrel s\over \to}\atop {\buildrel t\over \to}}\atop {\buildrel \hookleftarrow \over i}} C,\quad si=ti=\id.\]
In this case  we note that $M$ is a $C$-bicomodule by via 
\[ \Delta_L=(t\tens\id)\Delta,\quad \Delta_R=(\id\tens s)\Delta\]
and $i$ becomes a bicomodule map.
\item an associative nonunital algebra structure $\circ:M\tens M\to M$ such that $\Delta,\eps$ are algebra homs.
\item the `unity condition'
\[ \circ(\id\tens i)\Delta_R=\circ( i\tens\id)\Delta_L=\id\]
and `source and target conditions' on $M\square M$: 
\[  s\circ=\eps\tens s,\quad t\circ=t\tens\eps\]
\item an `antipode' $\CS:M\to M$ obeying the `twisted morphism' condition
 \[ (s\tens \CS)\Delta=\Delta_L\CS,\quad (\CS\tens t)\Delta=\Delta_R\CS\]
 and 
 \[  \circ(\CS\tens\id)\Delta=is,\quad \circ(\id\tens \CS)\Delta= i t\]
\end{enumerate}
\end{definition}

We call our groupoid-like object  `embedded' since the algebra product $\circ$ is only really wanted on $M\square M$.  However, at least in our example $\circ$ is defined on all of $M\tens M$ and this is convenient for the formulation of the homomorphism property of $\Delta$.  This is a working definition in that these are all properties that hold in our example but the full list of axioms is meant to emerge later as properties of such objects are required. For example, in our example $\circ$ on $M\tens M$ has a left unit $1_L$, but it is not clear that this is important so we have not explicitly included it in our list above.  

\begin{proposition}  Suppose  (1)-(2) in Definition~2.3. 

(i)  The conditions (3) are equivalent to requiring the restriction $\circ:M\square M\to M$  to be a morphism in the category of ${}^C\CM^C$ of bicomodules and form an associative algebra in this category and with unit $i$.

(ii) The conditions (3),(4) only refer to $\circ$ restricted to $M\square M$.

(iii) The conditions (3),(4) imply that
\[ \eps \CS=\eps, \quad t\CS=s,\quad s\CS=t\]
\end{proposition}
\proof (i) The requirement for $\circ$ to be morphism in ${}^C\CM^C$ when restricted to $\sum m\tens n\in M\square M$ is
\[ \sum t(m\o)\tens m\t\circ n=\sum t((m\circ n)\o)\tens (m\circ n)\t\]
\[ m\tens n\o\tens s(n\t)= \sum (m\circ n)\o\tens s((m\circ n)\t)\]
which by (2) means
\[ \sum t(m\o)\tens m\t\circ n=\sum t(m\o\circ n\o)\tens m\t \circ n\t\]
\[ \sum m\tens n\o\tens s(n\t)= \sum m\o \circ n\o\tens s(m\t \circ n\t)\]
These hold if the  $s\circ $ and $t\circ $ properties in axioms (3) hold, and conversely if they hold then apply $\eps$ to the right hand sides in the relevant factors
to conclude the $s\circ $ and $t\circ $ properties.

Next, suppose $\circ$ is a morphism and consider its restriction to $M\square M\to M$. SInce it is associative on $M\tens M$ by condition (2), it remains associative on this subspace and we have equality of the two iterated $\circ$ maps $M\square M\square M\to M\square M\to M$ (where the two brackettings are identified by abuse of notation, more precisely they are isomorphic by the underlying associativity of vector space $\tens$). Finally, note that $C\square M\isom M$ and  $M\square C\isom M$ by $\Delta_L,\Delta_R$ respectively (in the reverse direction use the counit of $M$) whereby $C$ is the unit object in the monoidal category with respect to $\square$. Hence  the requirement of  $\circ$ to be unital under $i:C\hookrightarrow M$ namely that $\circ(i\tens\id): C\square M\to M\square M\to M$  and $\circ(\id\tens i):M\square C\to M\square M\to M$ are the identity after allowing for these isomorphisms, becomes the `unity conditions' stated in (3). 

(ii) By the explanation given for $C\square M\isom M$ and  $M\square C\isom M$ and $i$ a bicomodule map, it is clear that $(\id\tens i)\Delta_R, (i\tens\id)\Delta_L$ have image in $M\square M$. Meanwhile, coassociativity and the twisted bicomodule property of $\CS$ ensures that the images of $(\CS\tens \id)\Delta, (\id\tens \CS)\Delta$ are also in $M\square M$:
\[  (\Delta_R\tens\id)(\CS\tens\id)\Delta =(\CS\tens t\tens\id)(\Delta\tens\id)\Delta=(\CS\tens t\tens\id)(\id\tens\Delta)\Delta=(\id\tens\Delta_L)(\CS\tens\id)\Delta\]
\[ (\Delta_R\tens\id)(\id\tens \CS)\Delta=(\id\tens s\tens \CS)(\Delta\tens\id)\Delta=(\id\tens s\tens \CS)(\id\tens\Delta)\Delta=(\id\tens\Delta_L)(\id\tens \CS)\Delta\]

(iii) For any embedded quantum groupoid we have $\eps \CS=\eps$ by applying $\eps$ to the second antipode axiom in (4). In this case applying $\eps$ to the twisted bicomodule condition we arrive at $t\CS=s, s\CS=t$. \endproof

 Next, if $C$ is a Hopf algebra and if 
 $M,N$ are $C$-bicomodule algebras (i.e. the coactions are algebra homs) then one may check that $M\square N$ gets the structure of a $C$-bicomodule algebra with the product $(m\tens n)(m'\tens n')=m m'\tens n n'$, i.e. the result lies in $M\square N$:
 \[ \Delta_R(mm')\tens nn'=\Delta_R(m)\Delta_R(m')\tens nn'=(\Delta_R(m)\tens n)(\Delta_R(m')\tens n')\]
 \[ =(m\tens\Delta_L n)(m'\tens \Delta_Ln')=mm'\tens \Delta_L(n)\Delta_L(n')=mm'\tens \Delta_L(nn').\]
 Then if $C,M$ are Hopf algebras and $s,t,i$ Hopf algebra maps in Definition~2.3 then $M$ becomes a $C$-bicomodule algebra under pushout via $s,t$ and hence $M\square M$ is a $C$-bicomodule algebra also.

 \begin{definition} We define a `strict quantum 2-group' as a pair of Hopf algebras and Hopf algebra maps 
\[ H_1 {{{\buildrel s\over \to}\atop {\buildrel t\over \to}}\atop {\buildrel \hookleftarrow \over i}} H_0,\quad  s i=\id,\quad t i=\id\]
together with an associative product  $\circ: H_1\tens H_1\to H_1$ forming an `embedded quantum groupoid' as in Definition~2.3 and subject to the `quantum interchange law'
\[  \circ \cdot_{H_1\square H_1}=\cdot_{H_1}(\circ\tens\circ)\]
as maps $(H_1\square H_1)\tens (H_1\square H_1)\to H_1$, where $\cdot_{H_1}$ is the Hopf algebra product of $H_1$ and $\cdot_{H_1\square H_1}$ is the product of the  cotensor $H_1\square H_1$ as bicomodule algebras. 
\end{definition}
Again, this is a working definition to lay out properties of our example rather than a fully refined list of axioms. In particular one might ask that the Hopf algebra unit of $H_1$ is a left-unit for $1_L$ as this will be the case in our example below. Note that the quantum groupoid antipode $\CS$ in Definition~2.3 should not be confused with the Hopf algebra antipode which also exists in the Hopf case.

\begin{theorem} Given a crossed module as in Definition~2.1 we have a `strict quantum 2-group' in the sense above with $H_1=A\lcross H$ and $H_0=H$. The product $\circ$ on $H_1$ and its relevant `antipode' are 
\[ (a\tens h)\circ (b\tens g)=\eps(h) ab\tens g,\quad \CS(a\tens h)=Sa\o\tens d(a\t)h,\quad \forall a,b\in A,\ h,g\in H\]
and $1_L=1\tens 1$ is a left unit for $\circ$.
\end{theorem}
\proof (i) We take $H_1=A\lcross H$ and $H_0=H$ and we check
\[ t((a\tens h)(b\tens g))=d(a h\o\la b) h\t g=d(a) h\o d(b) Sh\t h\th g=d(a)h d(b) g=t(a\tens h)t(b\tens g)\]
\[ \Delta t(a\tens h)=d(a)\o h\o \tens d(a)\t h\t=d(a\o) h\o\tens d(a\t) h\t=t(a\o\tens h\o)\tens t(a\t\tens h\t)\]
as required. The other maps $i(h)=1\tens h$ and $s(a\tens h)=\eps(a)h$ are already known to be Hopf algebra maps with $s i=\id$. We similarly have $t(1\tens h)=d(1)h=h$. 

(ii) We check that the proposed product $\circ$ is associative and left unital with $1_L=1\tens 1$:
\begin{eqnarray*}  ((a\tens h)\circ (b\tens g))\circ(c\tens f)&=&\eps(h)(ab\tens g)\circ(c\tens f)=\eps(h)\eps(g) abc\tens f\\
(a\tens h)\circ((b\tens g)\circ(c\tens f))&=&(a\tens h)\circ (bc\tens f)\eps(g)=abc\tens f\eps(h)\eps(g)\\
(1\tens 1)\circ (a\tens h)&=& a\tens h \eps(1),\quad (a\tens h)\circ(1\tens 1)=\eps(h)a\tens 1.\end{eqnarray*}
We also check that $\Delta$ is an algebra hom with respect to $\circ$ and the Hopf algebra coproduct of $A\lcross H$:
\begin{eqnarray*}  \Delta((a\tens h)\circ(b\tens g))&=&a\o b\o\tens g\o \tens a\t b\t\tens g\t\eps(h)\\ 
&=& a\o b\o\tens g\o\eps(h\o)\tens a\t b\t\tens g\t\eps(h\t)\\
&=&(a\o\tens h\o)\circ (b\o\tens g\o)\tens (a\t\tens h\t)\circ(b\t\tens g\t)\end{eqnarray*} 
and similarly for the counit. 

(iii) From the coactions, the space $H_1\square H_1$ is
\[ H_1\square H_1=\{\sum a\tens h\tens b\tens g\ |\ \sum a\tens h\o\tens h\t\tens b\tens g=\sum a\tens h\tens d(b\o)g\o\tens b\t\tens g\t\} \]
where the sum is to remind us that elements will typically be sums of the form shown. Applying $\eps$ to both sides we variously deduce in this case
\[ \sum a\tens h\tens b\tens g= \sum a\tens \eps(h) d(b\o)g\o\tens b\t\tens g\t\]
\[ \sum a\tens \eps(b)h\tens g=\sum a\tens \eps(h) d(b)g\o\tens g\t\]
\[ \sum a\tens \eps(b)\eps(g)h= \sum a\tens \eps(h) d(b)g\]
which will be needed below for elements of this space. We now verify that the restricted $\circ$ is a morphism in ${}^H\CM^H$. For left covariance:
\begin{eqnarray*}\sum t(a\o\tens h\o)\tens  (a\t \tens h\t)\circ (b\tens g)&=&\sum d(a\o)h\tens a\t b \tens g\\
\Delta_L(\sum(a\tens h)\circ(b\tens g))=\sum \Delta_L(ab\tens g\eps(h))&=&\sum d(a\o b\o)g\o\tens a\t b\t\tens g\t \eps(h)
\end{eqnarray*}
which are equal by multiplicativity of $d$ and the first of our identities above for elements of $H_1\square H_1$. For right covariance:
\begin{eqnarray*}&&\sum (a\tens h)\circ(b\o\tens g\o)\tens s(b\t\tens g\t)= \sum (a\tens h)\circ (b\tens g\o)\tens g\t\\
&&=\sum a b\tens g\o \tens g\t\eps(h)=\sum a\o b\o\tens g\o \eps(h)\tens s(a\t b\t\tens g\t)\\
&&=\sum\Delta_R(ab\tens g\eps(h))=\sum \Delta_R(\sum(a\tens h)\circ(b\tens g))\end{eqnarray*}
(in fact even without restriction). One can similarly verify that
\begin{eqnarray*} &&\sum t((a\tens h)\circ  (b \tens h)):=\sum t((ab\tens g)\eps(h))=\sum \eps(h)d(ab)g\\
&&=\sum d(a)\eps(h) d(b)g =\sum d(a)\eps(b)\eps(g)h=: \sum t(a\tens h)\eps(b\tens g)\end{eqnarray*}
where the fourth step is the last of the identities for an element of $H_1\square H_1$ in (iii), and
\begin{eqnarray*} &&\sum s((a\tens h)\circ  (b \tens h)):=\sum s((ab\tens g)\eps(h))=\sum \eps(h)\eps(a)\eps(b) g=: \sum \eps(a\tens h)s(b\tens g)\end{eqnarray*}
in fact on all of  $H_1\tens H_1$. These identities are implied, however, from the bicovariance of $\circ$ (apply the counit of $H_1$). 
Similarly, as $\circ$ on $H_1\tens H_1$ is associative it necessarily restricts to an associative algebra in ${}^{H_1}\CM{}^{H_1}$ on $H_1\square H_1$ in the sense that the two iterated $\circ$ maps $H_1\square H_1\square H_1\to H_1\square H_1\to H_1$ (using the vector space associator implicitly to identify the two brackettings of $\square$) coincide. Finally, we verify the unit map $i$: 
\[ \circ(\id\tens i)\Delta_R(a\tens h)=(a\o\tens h\o)\circ(1\tens s(a\t\tens h\t))=(a\tens h\o)\circ(1\tens h\t)=a\tens h\]
\[ \circ(i\tens\id)\Delta_L(a\tens h)=(1\tens t(a\o\tens h\o))\circ(a\t\tens h\t)=(1\tens d(a\o)h\o)\circ(a\t\tens h\t)=a\tens h\]
for all $a\tens h\in H_1$. 

(iv) We verify the twisted bicomodule properties 
\begin{eqnarray*} (s\tens \CS)\Delta(a\tens h)&=&s(a\o\tens h\o)\tens Sa\t\tens d(a\t)h\t=h\o\tens Sa\o\tens d(a\t)h\t\\
\Delta_L\CS(a\tens h)&=&\Delta_L(Sa\o\tens d(a\t)h)\\
&=&t(Sa\o\t\tens d(a\t\o)h\o\tens Sa\o\o\tens d(a\t\t)h\t\\
&=&d(Sa\o\t) d(a\t\o)h\o\tens Sa\o\o\tens d(a\t\t)h\t \end{eqnarray*}
which agree using the antipode axioms and $d$ and algebra map. On the other side
\begin{eqnarray*} (\CS\tens t)\Delta(a\tens h)&=&S(a\o\tens h\o)\tens d(a\t)h\t=Sa\o\tens d(a\t)h\o\tens d(a\th)h\t\\
\Delta_R\CS(a\tens h)&=&\Delta_R(Sa\o\tens d(a\t)h)=Sa\o\t\tens d(a\t\o)h\o\tens \eps(a\o\o)d(a\t\t)h\t
\end{eqnarray*}
which agree. We next verify  the antipode properties
\begin{eqnarray*} \CS(a\o\tens h\o)\circ (a\t\tens h\t)&=&(Sa\o \tens d(a\t)h\o)\circ (a\th\tens h\t)\\
&=&(Sa\o)a\t\tens h)= \eps(a)1\tens h=is(a\tens h)\\
(a\o\tens h\o)\circ \CS(a\t\tens h\t)&=&(a\o\tens h\o)\circ( Sa\t\tens d(a\th)h\t)\\
&=& a\o Sa\t\tens d(a\t)h=1\tens d(a)h=it(a\tens h)\end{eqnarray*}
One can also check that 
\[ t\CS(a\tens h)=t(Sa\o\tens d(a\t)h)=d(Sa\o)d(a\t)h=\eps(a)h=s(a\tens h)\]
as $d$ is an algebra map, which was mentioned as implied by the general theory.

(v) Finally, we verify the interchange law
\begin{eqnarray*}&&\kern-20pt \sum ((a\tens h)(a'\tens h'))\circ ((b\tens g)(b'\tens g'))\\
&&:=\sum (a(h\o\la a')\tens h\t h')\circ(b (g\o\la b')\tens g\t g')\\
&&=\sum a(h\la a')b(g\o\la b')\tens g\t g'\eps(h')\\
&&=\sum a(d(b\o)g\o\la a')b\t(g\t\la b')\tens g\th g'\eps(h)\eps(h')\\
&&=\sum a b\o (g\o\la a')(Sb\t)b\th(g\t\la b')\tens g\th g'\\
&&=\sum ab(g\o\la a')(g\t\la b')\tens g\th g'\eps(h)\eps(h')\\
&&=\sum (ab\tens g)(a' b'\tens g')\eps(h)\eps(h')\\
&&=:\sum ((a\tens h)\circ (b\tens g))((a'\tens h')\circ (b'\tens g'))\end{eqnarray*}
where the third equality uses $\sum a\tens h\tens b\tens g\in H_1\tens H_1$ and the first of the displayed identities from this in (iii). The fourth
equality is our condition (3) of Definition~2.1. Elsewhere we use the product of $H_1=A\lcross H$ and the $\circ$ product. The proof  holds more generally
on any elements $\sum a'\tens h'\tens b'\tens g'\in H_1\tens H_1$. 
\endproof

This theorem exactly includes the classical case covered in Section~1 when the Hopf algebras are taken to be group algebras but also suggests the possibility of a still more general construction in the next section. We also observe in this regard,

\begin{proposition} For the strict quantum 2-groupoids in Theorem~2.6 we have:

(i)  $\sum \CS(m\circ n)=\sum \CS(n)\circ \CS(m)$ for all $\sum m\tens n\in H_1\square H_1$ if 
\[ a\o\tens d(a\t)= a\t\tens d(a\o),\quad \forall a\in A\]

(ii)  $\CS^2=\id$ iff the condition displayed in (i) and $A$ is involutive.

(iii) $\tau(\CS\tens \CS)\Delta= \Delta \CS$ where $\tau$ is the flip map iff the condition in (i) and $H$ is cocommutative.

(iv) $A\tens 1, 1\tens A^H\subset H_1\square H_1$ where ${}^HA=\{a\in A\ |\ d(a\o)\tens a\t=1\tens a\}$. 
\end{proposition}
\proof (i) We compute
\begin{eqnarray*} \CS((a\tens h)\circ b\tens g))&=&\CS(ab\tens g\eps(h))=S(a\o b\o)\tens d(a\t)d( b\t)g\eps(h)\\
(\CS(a\tens h))\circ \CS(b\tens g)&=& (Sb\o\tens d(b\t)g)\circ (Sa\o\tens d(a\t)h)=S(a\o b)\tens d(a\t)h\eps(g)\end{eqnarray*}
Now if $\sum a\tens h\tens b\tens g\in H_1\square H_1$ then by applying $\eps$ to the first of the  identities for this in the proof of Theorem~2.6, we have
\[ \sum a\tens h\tens b\eps(g)=\sum a\tens d(b\o)g\tens b\t \eps(h)\]
and hence if the condition displayed in part (i) holds then
\[ \sum a\tens b\tens h\eps(g)=\sum a\tens b\o\tens d(b\t)g\eps(h)\]
and we see that $S$ is then anti-multiplicative on such elements. 

(ii) Using $d$ an algebra hom and the antipode axioms on $A$ we find easily that
\[ \CS^2(a\tens h)=S^2(a\t)\tens d(Sa\o a\th)h,\quad \forall a\in A,\ h\in H\]
Hence if the condition displayed in (i) holds then $\CS^2(a\tens h)=S^2(a)h$  (use the condition in the form $a\t\tens d(a\th)=a\th\tens d(a\t)$ then
cancel the antipode. Conversely, if $\CS^2=\id$ then applying to $h=1$ we have $S^2(a\t)\tens d(Sa\o a\th)=a\tens 1$ for all $a\in A$. By applying the
counit we see that $A$ is involutive and in this case $a\tens 1=a\t\tens d(Sa\o a\th)$ for $a\in A$. Applying this to $a\o$ in the expression $a\t \tens d(a\o)1$ we have
$a\t\tens d(a\o)=a\th\tens d(a\o)d(Sa\t a_{(4)})=a\o\tens d(a\t)$ using $d$ an algebra map and cancelling the antipode. 

(iii) We compute, using that $d$ is a coalgebra map,
\begin{eqnarray*}\Delta \CS(a\tens h)&=&\Delta (Sa\o\tens d(a\t)h)\\
&=&Sa\o\t\tens d(a\t)\o h\o\tens Sa\o\o\tens d(a\t)\t h\t\\
&=&Sa\t\tens d(a\th)h\o\tens Sa\o\tens d(a_{(4)})h\t\\
\tau(\CS\tens \CS)\Delta(a\tens h)&=&\CS(a\t\tens h\t)\tens \CS(a\o\tens h\o)\\
&=&Sa\t\o\tens d(a\t\t)h\t\tens Sa\o\o\tens d(a\o\t)h\o\\
&=&Sa\th\tens d(a_{(4)})h\t\tens Sa\o\tens d(a\t) h\o\end{eqnarray*}
Comparing we see that
two applications of the condition stated in (i) and $H$ cocommutative imply equality. Conversely, if equality then setting $a=1$ implies $H$ cocommutative
and setting $h=1$ and applying the counit in two places gives the condition stated in (i) (assuming $A$ has invertible antipode). 

(iv) This is immediate from the characterisation of $H_1\square H_1$ in the proof of Theorem~2.6 to elements of the form $a\tens 1\tens 1\tens 1$ and $1\tens 1\tens a\tens 1$.  Note that ${}^HA\subset A^{inv}$ by applying the counit (we also know that $(i\tens\id)\Delta_L(H_1), (\id\tens i)\Delta_R(H_1),(\CS\tens\id)\Delta(H_1), (\id\tens \CS)\Delta(H_1)\subset H_1\square H_1$ from Proposition~2.4 (ii)).   \endproof

Finally, for ease of comparison with the more familiar notion of bialgebroid or Hopf algebroid\cite{Tak,Sch,Bo,BS}, we conclude here with the dual notion to Definition~2.3. Thus we replace $C$ by an algebra $R$ and $M$ by an algebra $A$. Reversing all arrows, we need algebra maps $s,t,\pi$ 
\[  R {{{\buildrel s\over \hookrightarrow}\atop {\buildrel t\over \hookrightarrow}}\atop {\buildrel \leftarrow \over \pi}} A,\quad \pi s=\pi t=\id.\]
The maps $s,t$ make $A$ an $R$-bimodule by $r.a.r'=t(r)as(r')$ for all $r,r'\in R$ and $a\in A$, and $\pi$ a bimodule map. Next, we need 
a non-counital coalgebra $\Delta:A\to A\tens A$ such that these maps are algebra homs.  Then we have a further projection to $\bar\Delta:A\to A\tens_R A$ and we require that this is a bimodule map and that $\pi$ is a counit for it. The former requirement is equivalent (given multiplicativity of $\Delta$ and hence of $\bar\Delta$) to
\[ \bar\Delta s(r)=1\tens_R s(r),\quad \bar\Delta t(r)=t(r)\tens_R1,\quad\forall r\in R\]
and the latter to 
\[ a\o s\pi(a\t)=t \pi(a\o)a\t=a,\quad\forall a\in A\]
where only the projected $\bar\Delta a=a\o\tens_Ra\t$ in view of the bimodule property of $\bar\Delta$. 

We then have $(\bar\Delta,\pi)$ a coalgebra in the category ${}_R\CM_R$. Finally we require an `antipode' $\CS:A\to A$ such that $\CS$ is a twisted bimodule map in the sense $\CS(t(r)a)=s(r)\CS a$, $\CS(as(r))=(\CS a)t(r)$ (this implies $\CS t=s$, $\CS s=t$ by setting $a=1$) and
\[ (\CS a\o)a\t=s\pi(a),\quad a\o \CS a\t=t\pi(a),\quad \forall a\in A\]
where only the projected $\bar\Delta a=a\o\tens_Ra\t$ is being used in these expressions in view of the twisted bimodule property and from where $\CS1=1$. Again this is a working definition, for example in the dual version of our main theorem we would also have a left counit $\eps_L$ for $\Delta$.  For the Hopf case $R,A$ will be Hopf algebras, $r,s,t$ will be Hopf algebra maps and there will be an interchange law involving the two coproducts on $A$. The key part of this is our observation that if $A,B$ are $R$-bimodules and comodule coalgebras under $R$ then so is $A\tens_R A$ with the tensor product coalgebra.

We now turn to the construction of examples. The first one is the obvious application:

\begin{example} There is an obvious notion of a Lie-crossed module obtained by differentiating the data in Section~1. This is given by a pair of Lie algebras  and a Lie algebra hom $d:\CA\to \CH$, a Lie algebra  action of $\CH$ on $\CA$, subject to 
\[ d(h\la a)=[h,d(a)],\quad d(a)\la b=[a,b],\quad\forall h\in \CH,\ a\in \CA.\]
We have a semidirect sum algebra $\CA\lcross\CH$ and Lie algebra maps $\CA\lcross\CH {{{\buildrel s\over \to}\atop {\buildrel t\over \to}}\atop {\buildrel \hookleftarrow \over i}}\CH$ given by 
\[ s(a)=0,\quad s(h)=h,\quad t(a)=d(a),\quad t(h)=h,\quad i(h)=h.\]
It is straightforward to exponentiate this data to $H=U(\CH)$ and $A=U(\CA)$ with actions and maps define as above on the generators. Hence we have a Hopf crossed module
and an easy `classical' example of a strict quantum 2-group by Theorem~2.6.  We have $H_1=U(\CA\lcross\CH)$ containing $U(\CH)$ (by $i$) and $U(\CA)$ as subHopf algebras. We denote the enveloping algebra products as usual by omission. The map $s$ sets the $\CA$ to zero and the map $t$ sends $a\in\CA$ to $d(a)\in \CH$, both leaving the $\CH$ generators unchanged. We have 
\[ a\circ X=aX,\quad h\circ X=0,\quad \forall a\in\CA,\ h\in\CH,\  X\in U(\CA\lcross\CH)\]
and the same for any monomial in $\CH$ of positive degree. Accordingly we cannot have a a right unit for $\circ$ but $1$ is a left unit.  Clearly Proposition~2.7 applies also in this case so $S$ is anti-co/multiplicative and involutive. We have $Sa=-a$, $Sh=h$ on generators. The coactions are 
\[ \Delta_La=d(a)\tens 1+1\tens a,\quad \Delta_Lh=h\tens 1+1\tens h,\quad \Delta_Ra=a\tens 1,\quad \Delta_Rh=h\tens 1+1\tens h\]
on generators. The subalgebra $H_1\square H_1$ depends on $d$ but is nontrivial as it contains, for example $\CA\tens 1$ and $1\tens \ker\, d$ from Proposition~2.7 (iv).
\end{example}

Clearly the same works for any pair of cocommutative Hopf algebras $A,H$ where $H$ acts on $A$ respecting its Hopf algebra structure, $d:A\to H$ a Hopf algebra map between them such that (2)-(3) in Definition~2.1 hold. 

\begin{example} Let $A=H$ a cocommutative Hopf algebra, $H$ acts on $A$ by the adjoint action and $d=\id$. This is easily seen to form a Hopf crossed module. We have a 2-quantum group with $H_0=H$ and $H_1=H\lcross H\isom H\tens H$ where the cross product by the adjoint action is known to be isomorphic to a tensor product\cite{Ma:book}. The map $H\lcross H\to H\tens H$ is $h \tens g\mapsto h g\o\tens g\t$ with inverse $h\tens g\mapsto h Sg\o\tens g\t$. Using this we compute the maps on the second form $H_1=H\tens H$ as a Hopf algebra with
\[ s(h\tens g)=\eps(h)g,\quad t(h\tens g)=h\eps(g),\quad i=\Delta\]
\[  (h\tens g)\circ (h'\tens g')=h(Sg)h'\tens g',\quad \CS(h\tens g)=g\tens h.\]
Meanwhile, the coations on $H_1$ in this form are $\Delta_L=\Delta\tens\id$, $\Delta_R=\id\tens\Delta$ from which we compute $H_1\square H_1=H\tens\Delta(H)\tens H$. Indeed, if $\sum a_i\tens h_i\tens b_i\tens g_i\in H_1\tens H_1$ we require 
\[ \sum a_i\tens h_i\o\tens h_i\t\tens b_i\tens g_i=\sum a_i\tens h_i\tens b_i\o\tens b_i\t\tens g_i\]
which on applying $\eps$ in the second position tells us that $\sum a_i\tens h_i\tens b_i\tens g_i=a_i\tens\Delta(\eps(h_i)b_i)\tens g_i$. Now, when restricted to $H_1\square H_1$ we have $(a\tens h\o)\circ(h\t\tens g)=a(Sh\o)h\t\tens g=\eps(h)a\tens g$. From this one can easily check the 2-quantum group interchange law as must be the case from the general theory.
\end{example}

Thus the adjoint action gives a relatively trivial tensor product 2-quantum group, as it should. This is the quantum version of the adjoint crossed module on a group leading after an isomorphism to $G_1=G\times G$ a trivial 2-group over $G_0=G$ where $i$ is the diagonal map and $s,t$ are the projection to 1st and 2nd factors respectively. The second product here is $(h,g)\circ(h',g')= (h,g')$ whenever $g=h'$ but we see that it is in our approach embedded in an associative product on all of $G_1\times G_1$ given (reading off from the above example) by $(h,g)\circ(h',g')=(hg^{-1}h',g')$.  Similarly for general group crossed modules: although we use Hopf algebra methods, when we apply Theorem~2.6 to group algebras we recover the usual theory as recalled in Section~1 with the main difference that $\circ$ is actually defined associatively on all of $G_1\times G_1$. 

For a purely `quantum' class of examples we can take $A$ a commutative Hopf algebra acted upon by $H$ a cocommutative one, as a module algebra and module coalgebra.  The only remaining condition for a crossed module of this type is the Hopf algebra map $d:A\to H$ intertwining the given action on $A$ with the adjoint one (condition (2) in Definition~2.1).  The 2-quantum group structure is read off from Theorem~2.6.

\begin{example} We take $H_0=H$ cocommutative, trivial $d=1\eps$ and $A=H^*$ in the finite dimensional case or appropriately dually paired in the infinite-dimensional case such that the coadoint action makes sense. There is a canonical coadjoint action of $H$ on $A$ and  $H_1=A\lcross H=D(H)$ is the quantum double Hopf algebra of $H$. In this class of examples $t=s$ and
\[ s(a\tens h)=\eps(a)h,\quad i(h)=1\tens h,\quad (a\tens h)\circ(b\tens g)=\eps(h)ab\tens g,\quad \CS(a\tens h)=Sa\tens h\]
according to Theorem~2.6. 

To see what this looks like explicitly, we can take  a finite group $G$, $A=k(G)$, $H=kG$ with action by ${\rm Ad}$ and trivial $d$. Then $H_1=D(G)$ the quantum double,  $t=s$ is the counit projection and 
\[ (a\tens h)\circ (b\tens g)=ab\tens g,\quad \forall h,g\in G,\ a,b\in k(G)\]
We can think of elements $A$ of $D(G)$ as $kG$-valued functions on $G$, i.e. $A=\sum_{g\in G}a_g g$ where $a_g\in k(G)$. Then for two such, $(A\circ B)=(\sum_{g\in G}a_g)B$ with pointwise product on the right. 

Similarly, we can take $A=k(SU_2)$, $H=U(su_2)$ acting by $\rm ad$ and trivial $d$. Then $H_1=D(U(su_2))$, $t=s$ and
\[ (a\tens h)\circ (b\tens g)=0,\quad \forall h\in su_2,\quad a,b\in k(SU_2),\ g\in U(su_2).\]
One can think of $H_1$ as generated by products in $U(su_2)$ of  $su_2$-valued functions on $SU_2$ and functions on $SU_2$, and if $A$ is an example of the former $B$ anything then $A\circ B=0$. If $A$ is an example of the latter and $B$ anything then $A\circ B=AB$ the pointwise product. \end{example}

And finally one of this `quantum' type i.e. with noncommutative and noncocommutative $H_1$ and with $d$ nontrivial:

\begin{example}  As data we take $G$ a locally compact Abelian group acting from the right by group automorphisms on a group $M$ and $\hat d:\hat G\to M^G$ a group hom from the Pontryagin dual group of $G$ to the fixed subgroup in $M$.  Taking a convolution algebra approach to $k G$, $f\in k(M)$ pulled back along $\hat d$ Fourier transforms to an element of $k G$ and this defines the map $d:k(M)\to k G$ for a suitable class of functions $k(M)$ (for example over $\C$). Here we stick to an algebraic treatment focussing on $\hat G$ and a commutative $\hat G$-graded Hopf algebra $A$ in the role of $k(M)$ and for simplicity we explain the case where $\hat G$ is finite, and denote it additively. Note that the grading is the same thing as an action of $k(\hat G)\isom k G$ via $(\delta_u f)(s)=\sum_{g\in G} u(g) f(s\ra g)$ for $s\in M$ and $\delta_u$ the Kronecker delta function at $u\in\hat G$, but we do not need this as we refer to the grading itself. Thus  $A=\oplus_{u\in\hat G}A_u$ and we write $f_u$ for the component of $f\in A$ in $A_u$. We have $H_0=k(\hat G)$ and $H_1=A\lcross k(\hat G)$ with product
\[ (f\tens\delta_u)(h\tens \delta_v)=f h_{u-v}\tens\delta_v,\quad \forall f,h\in A,\ u,v\in
\hat G\]
The rest of the structure in Theorem~2.6 is, with $e\in M$ the group identity, 
\[ s(f\tens \delta_u)=f(e)\delta_u,\quad t(f\tens \delta_u)=f(\hat d(u))\delta_u\]
\[ (f\tens \delta_u)\circ(h\tens\delta_v)=\delta_{u,0}fh\tens\delta_v,\quad \CS(f\tens\delta_u)=L_{\hat d(u)^{-1}}(S f)\tens \delta_{u}\]
where $L_s(f)=f(s(\ ))$ is left translation by $s\in M$ and $(Sf)(s)=f(s^{-1})$ is given by inversion on $M$. 
The left and right coactions of $k(\hat G)$ mean respectively right and left actions of $\hat G$ and these are
\[  (f\tens\delta_v)\ra u=L_{\hat d(u)}(f)\tens\delta_{v-u},\quad u\la(f\tens \delta_v)=f\tens\delta_{v-u}.\]
If we also define $\psi:k(\hat G)\tens A\to A\tens k(\hat G)$ by
\[ \psi(f\tens\delta_u)=\delta_u\tens L_{\hat d(u)}(f),\quad\forall u\in \hat G,\ f\in A\]
then one may determine after some computation that 
\[ H_1\square H_1=(\id\tens \psi\tens\id)(A\tens A\tens \Delta(H))\]
in this class of examples. One may also check that $\CS^2(f\tens\delta_u)=f(\hat d(u)^{-1}(\ )\hat d(u))\tens\delta_u$ so $\CS^2=\id$ precisely when $\hat d(\hat G)$ is central in $M$, in agreement with the content of Proposition~2.7. 

To give a concrete but infinite example one can take $A=\C[SU_2]$ with its usual generators $a,b,c,d$ for the matrix entry functions (subject to $ad-bc=1$). We take $\hat G=\Z$ and the grading that corresponds to $S^1\subset SU_2$ acting from the right by conjugation, namely $|a|=|d|=0$, $|b|=-2$, $|c|=2$. The fixed subgroup in the unitary setting is $S^1$ and we identify  $\omega$ a complex number of modulus $1$ with ${\rm diag}(\omega,\omega^{-1})$ in $S^1$. So the map $\hat d$ means to specify an arbitrary $\omega\in S^1$ so that $\hat d(n)=\omega^n\in S^1$ viewed in $SU_2$. Then $t(f\tens\delta_n)=f(\omega^n)\delta_n$ compared with $s(f\tens \delta_n)=f(e)\delta_n$. Finally, $A$ has another grading $\deg(a)=\deg(b)=1$, $\deg(c)=\deg(d)=-1$ familiar from the Hopf fibration so that $L_{\hat d(n)}(f)=\omega^{n\deg(f)}f$ on homogeneous elements. From these facts one can determine that $\CS$ has finite order iff $\omega$ is a root of unity and in this case the order is even. Indeed, $\CS^{2n}=\id$ if $\omega^{2n}=1$ in view of the formula above for $\CS^2$. Note that $k(\Z)$ is not a Hopf algebra due to the infinite sum in its coproduct but most of the other structure still makes sense by working with actions as we have done above (there are standard methods to handle the coproduct as well, notably the theory of locally compact quantum groups).  Clearly, this example works similarly for all compact Lie groups $M$ with $G\subseteq M$ the maximal torus and $\hat d$ defined by ${\rm rank}(M)$ elements of $G$. \end{example}

\section{Braided Hopf crossed modules}

The characterization of Proposition~2.2 suggests the following more general definition which we study here, although we do not have a full quantum groupoid picture in this more general case. We assume $H$ has invertible antipode. A Hopf algebra $B\in\CZ({}_H\CM)$ then means first of all an algebra and coalgebra where both structure maps are morphisms i.e. equivariant under the action and coaction of $H$. We refer to the latter as the `given' action and coaction and we denote them by $\la$ and $b\bo\tens b\bt$ respectively.  In addition the coproduct is a homomorphism which translates as
\[ \Delta (ab)=a\o\Psi(a\t\tens b\o)b\t=a\o (a\t\bo\la b\o)\tens a\t\bt b\t,\quad\forall a,b\in B\]
and there is also a braided antipode $S:B\to B$. We refer to \cite{Ma:book} for the general theory of Hopf algebras and \cite{Ma:bg,Ma:introm} for braided-Hopf algebras or Hopf algebras in braided categories. The trivial given coaction would take us back to $B$ an ordinary Hopf algebra and the theory of Section~2.

\begin{definition} A braided Hopf crossed module means
\begin{enumerate}
\item A Hopf algebra $H$ and a braided-Hopf algebra  $B$ in the category $\CZ({}_H\CM)$
\item A `twisted Hopf algebra map'  $d:B\to H$ meaning an algebra map obeying
\[  \Delta d(b)= d(b\o)b\t{}\bo\tens d(b\t\bt),\quad \eps(d(b))=\eps(b),\quad \forall b\in B\]
\item $d$ also obeys 
\[ d(h\la b)= h\o d(b) Sh\t,\quad \forall b\in B,\ h\in H.\]
\item 
\[ d(a)\la b=a\o (a\t \bo\la b) Sa\t\bt,\quad\forall a,b\in B\]
(the right hand side is the braided adjoint action in the category $\CZ({}_H\CM)$). Without condition (4) we say that we have a braided pre-crossed module. 
\end{enumerate}
\end{definition}

The data (1) means that we have a biproduct ordinary Hopf algebra
\[  B\lbiprod H {{\buildrel s\over \to}\atop {\buildrel \hookleftarrow \over i}} H\]
where we make the smash product by the given action of $H$ and also a smash coproduct by the
given coaction, i.e.
\[ \Delta(b\tens h)= b\o \tens b\t\bo h\o \tens b\t\bt\tens h\t\]
and the maps shown are the obvious ones from the unit and counit. We similarly have an inclusion
\[ B\hookrightarrow B\lbiprod H\]
but only as algebras. By a theorem that goes back to Radford\cite{Rad}, this data is the most general case for a larger Hopf algebra with projection (in the sense of Hopf algebra inclusion and covering surjection) to a smaller one. The braided-category interpretation of this theorem appeared only after the invention of braided-Hopf algebras and is in \cite{Ma:skl}.  This is another motivation for a generalisation of Section~2, namely a Hopf algebra with projection as in Definition~2.5 necessarily has the form $B\lbiprod H$ for some $B\in \CZ({}_H\CM)$. Applying this in Theorem~2.6 with projection via $s$ gives back the ordinary Hopf algebra $A$ but in general the converse should be associated with a braided one. The other projection via $t$ should also play a role, to be explored elsewhere.

\begin{lemma} (i) Given (1) the data (2) in Definition~3.1 allow us to define an `induced coaction' on $B$
\[ \Delta_L:B\to B\tens H,\quad \Delta_L(b)=d(b\o)b\t\bo\tens b\t\bt\]
(ii) The  given action $\la$ and $\Delta_L$ fit together to make $B$ an object of $\CZ({}_H\CM)$  iff the
condition (3) in Definition~3.1 holds. (iii) Under these conditions $B$ is an algebra in $\CZ({}_H\CM)$ and 
$d:B\to H$ is a morphism.
\end{lemma}
\proof (i) We check first that $\Delta_L$ as stated is a coaction:
\begin{eqnarray*}  (\id\tens\Delta_L)\Delta_Lb &:=& d(b\o)b\t\bo\tens d(b\t\bt\o)b\t\bt\t\tens b\t\bt\t\bt\\
&=&d(b\o)b\t\o\bo b\t\t\bo\tens \tens d(b\t\o\bt)b\t\t\bt\bo\tens b\t\t\bt\bt\\
&=&d(b\o)b\t\o\bo b\t\t\bo\o\tens \tens d(b\t\o\bt)b\t\t\bo\t \tens b\t\t\bt\\
&=& d(b\o)\o b\t\bo\o\tens  d(b\o)\t b\t\bo\t\tens b\t\bt=:(\Delta\tens\id)\Delta_Lb\\
 \end{eqnarray*}
 using in the second equality  that the coproduct $\Delta$ of $B$ is equivariant under the given (original) coaction of $H$, meaning
an $H$-comodule coalgebra
\[ b\bo\tens\Delta b\bt= b\o\bo b\t\bo \tens b\o\bt\tens b\t\bt.\]
The third equality is that the given coaction is indeed a coaction. The fourth is the twisted Hopf algebra map assumption for $d$ and after that we recognise the answer.  

(ii) Next assuming (3) in Definition~3.1 we check that we have an crossed$H$-module with $\Delta_L$ and the given action on $B$,
\begin{eqnarray*} \Delta_L(h\la b)&:=&d((h\la b)\o) (h\la b)\t\bo\tens (h\la b)\t\bt \\
&=& d(h\o\la b\o)(h\t\la b\t)\bo\tens (h\t\la b\t)\bt\\
&=& h\o\o d(b\o)(Sh\o\t) h\t\o b\t\bo Sh\t\th \tens h\t\t\la b\t\bt\\
&=& h\o d(b\o) b\t\bo Sh\th \tens h\t\la b\t\bt  \end{eqnarray*}
which is the required right hand side of the crossed $H$-module condition for coaction $\Delta_L$. The second equality is that the coproduct of $B$ is equivariant under the given action of $H$. The third equality is the condition (3) of Definition~3.1 after which we cancel the antipode. Conversely the holding of the 3rd equality implies the condition (3) of Definition~3.1 by applying the counit to the second factors and cancelling the antipode. 

(iii) Under these conditions we check that $\Delta_L$ is an algebra hom 
\begin{eqnarray*} \Delta_L(ab)&=& d(a\o (a\t\bo\la b\o))(a\t\bt b\t)\bo\tens (a\t\bt b\t)\bt\\
&=& d(a\o)a\t\bo\o d(b\o)(Sa\t\bo\t)a\t\bt\bo b\t\bo\tens a\t\bt\bt b\t\bt\\
&=& d(a\o)a\t\bo\o d(b\o)(Sa\t\bo\t)a\t\bo\th b\t\bo\tens a\t\bt b\t\bt\\
&=& \Delta_L(a)\Delta_L(b) \end{eqnarray*}
where the first equality uses the braided coproduct of $ab$ in computing $\Delta_L(ab)$. The second equality the condition (3) of Definition~3.1 already used as well as that $B$ is a comodule algebra for the original given coaction.  The third equality is the coaction property of the latter and coassociativity to renumber. We then cancel the antipode and recognise the answer. Hence $B\in \CZ({}_H\CM)$ as an algebra with respect to the new crossed $H$-module structure afforded by $\Delta_L$. 

Finally, we verify that also in the above case $d$ becomes a morphism in $\CZ({}_H\CM)$ with respect to this new induced coaction. The action is unchanged and equivariance here is just (2) in Definition~3.1. For the coaction we need
\[ (\id\tens d)\Delta_L(b):=d(b\o)b\t\bo\tens d(b\t\bt)=\Delta d(a)\]
which is just condition (1) in Definition~3.1. Thus we can deduce the `twisted Hopf algebra map' requirement from this morphism property if one assumes the form of $\Delta_L$ (or vice versa as we have done). \endproof

In the case of the trivial map $d=1\eps$  the induced coaction $\Delta_L$ coincides with the original given coaction so in this case the lemma does not produce anything new, but for general $d$ we end up with a different coaction of $H$ and when both parts hold we end up with $B\in \CZ({}_H\CM)$ as an algebra in a different way than the original given one. 

\begin{lemma} (i) Given a braided precrossed module, the remaining condition (4) in Definition~3.1 is equivalent to $B$ braided commutative as an algebra in $\CZ({}_H\CM)$ when viewed with the induced coaction $\Delta_L$. (ii) In this case the  subalgebra $B^{inv}=\{b\in B\ |\ d(b)=\eps(b)\}$ is commutative in $\CZ({}_H\CM)$ when viewed with the original given coaction.
\end{lemma}
\proof (i) The braiding in $\CZ({}_H\CM)$ with the new induced coaction $\Delta_L$ in Lemma~3.2 is
\[ \Psi_{\Delta_L}(a\tens b)=(d(a\o)a\t\bo)\la b\tens a\t\bt\]
which if condition (4) holds becomes
\[ \Psi_{\Delta_L}(a\tens b)=a\o(a\t\bo\la b)Sa\t\bt\o\tens a\t\bt\t.\]
Hence in this case we find
\[ \cdot\Psi_{\Delta_L}(a\tens b)=a\o(a\t\bo\la b)(Sa\t\bt\o) a\t\bt\t=ab\]
on cancelling the antipode, so the algebra is braided-commutative. Conversely, if
\[ ab=((d(a\o)a\t\bo)\la b)a\t\bt\]
then
\begin{eqnarray*}(a\o(a\t\bo\la b))Sa\t\bt&=&(d(a\o\o)a\o\t\bo a\t\bo)\la b a\o\t\bt Sa\t\bt\\
&=&(d(a\o)a\t\o\bo a\t\t\bo )\la b\, a\t\o\bt S a\t\t\bt\\
&=&( d(a\o)a\t\bo\la b)\, a\t\bt\o S a\t\bt\t=d(a)\la b\end{eqnarray*}
on cancelling the antipode. Hence in this case the condition (4) holds. 

(ii) Clearly 
\[ \Psi_{\Delta_L}(a\tens b)=a\bo\la b\tens a\bt,\quad\forall a\in B^{inv},\quad b\in B\]
so for elements of $B^{inv}$ the braiding reduces to that for the original given structure of $B\in \CZ({}_H\CM)$. Hence when condition (4) in Definition~3.1 holds we see that $B^{inv}$ is a commutative in the original braided category where $B$ was given as a braided Hopf algebra. \endproof

\begin{proposition} A braided Hopf crossed module is equivalent to the following data
\begin{enumerate}
\item $B$ a Hopf algebra in the category $\CZ({}_H\CM)$ over a Hopf algebra $H$.
\item  $d:B\to H$ is a `twisted Hopf algebra map' in the sense of (2) in Definition~3.1 such that the given action and the induced coaction make $B$ an object of  $\CZ({}_H\CM)$ (as such $B$ is an algebra and $d$ becomes a morphism).
\item $B$ is commutative when viewed in $\CZ({}_H\CM)$ according to the induced coaction. \end{enumerate}
\end{proposition}
\proof This brings together the above results. The first item is unchanged from Definition~3.1. The second item is (ii) of Lemma~3.2 where we continue to assume (i). There are other formulations as mentioned in the proof of the lemma. The third item is (i) in Lemma~3.3. \endproof

What we would have been nice at this point is a theorem that a braided crossed module leads to a strict quantum 2-group $H_1 {{{\buildrel s\over \to}\atop {\buildrel t\over \to}}\atop {\buildrel \hookleftarrow \over i}} H_0$ in the sense of Definition~2.5 with $H_1=B\lcross H$ and $H_0=H$. In this greater generality we appear to  have only part of this structure:

\begin{proposition} Given a braided Hopf crossed module $B$ in $\CZ({}_H\CM)$ we have a Hopf algebra $B\lbiprod H$ and a  projection given by $s(b\tens h)=\eps(b)h$ and $i(h)=1\tens h$. In addition $t(b\tens h)=d(b)h$  is another Hopf algebra map giving a projection so that $ti=\id$. 
\end{proposition}
\proof The Hopf algebra is the Radford biproduct in our braided-interpretation \cite{Ma:skl}, the new part is the map $t$. Of this, the algebra homomorphism part is unchanged from the proof of Theorem~2.6 so we have only to check the coalgebra part:
\begin{eqnarray*}
\Delta t(a\tens h)&:=& d(a)\o h\o \tens d(a)\t h\t= d(a\o)a\t\bo h\o \tens d(a\t\bt) h\t\\
&=& (t\tens t)(a\o\tens a\t\bo h\o\tens a\t\bt\tens h\t)=:(t\tens t)\Delta (a\tens h)\end{eqnarray*}
where we used the twisted-Hopf algebra condition in (2) of Definition~3.1 for the second equality and recognise
$(t\tens t)$ of the cross product coalgebra of $B\lbiprod H$. \endproof

It does not appear to be straightforward to find a second product $\circ$ and `antipode' so as to exactly fit  Definition~2.5 and it seems likely that the latter will need to be generalised in the braided case. This will be explored elsewhere. We give an example of a braided crossed module.

\begin{example} Let $H$ be a quasitriangular Hopf algebra so that its category of modules is braided, and $B\in {}_H\CM$ a braided Hopf algebra. In this case the functor \cite{Ma:dou} from ${}_H\CM\to \CZ({}_H\CM)$ allows us to view $B$ as required, with  left coaction $b\bo\tens b\bt=\CR\bt\tens \CR\bo\la b$ where $\CR=\CR\bo\tens\CR\bt$ is the quasitriangular structure. In this case $B\lbiprod H$ is the standard bosonisation of $B$ in \cite{Ma:bos}. In particular, we can take $B=B(H)$ the braided version of $H$ by transmutation. This has the same algebra as $H$ but a new coproduct $\und\Delta h=h\o S\CR\bt\tens \CR\bo\la h\t$ in terms of the initial coproduct of $H$. Next we take $d=\id$
and check the right hand side of the condition in (2) of Definition~3.1
\[ d(b\o)b\t{}\bo\tens d(b\t\bt)=h\o S\CR\bt (\CR\bo\la h\t)\bo\tens (\CR\bo\la h\t)\bt\]
\[=h\o S\CR\bt \CR'\bt \tens \CR'\bo\CR\bo\la h\t=h\o\tens h\t\]
since $(S\tens\id)\CR=\CR^{-1}$ and $(S\tens S)\CR=\CR$. The condition (3) in Definition~3.1 also holds since the given action of $H$ on $B$ is by the quantum adjoint action. Hence we have a braided precrossed module. We also see by the above computation that $\Delta_L=\Delta_H$ is the  coaction induced by $d=\id$ in Lemma~3.2 and it follows that $B,$ which has the same algebra as $H$, is commutative with respect to the associated braiding given  by this and the adjoint action. Hence by Lemma~3.3 we have a braided crossed module. Note that in this case $B\lbiprod H\isom H\codcross H$ as explained in \cite[Thm 7.4.6]{Ma:book} where the right hand side has tensor product algebra and the twisted coproduct $\CR_{32}\Delta_{H\tens H}(\ )\CR_{32}^{-1}$. 
In this form of $H_1=H\codcross H$ the maps $s,t,i$ are the same as in Example~2.9. The coactions $\Delta_{L,R}$ induced by $t,s$ in Defintion~2.3 also work out the same and hence $H_1\square H_1$ is the same as a vector space as in Example~2.9. The same $\circ$, however, is not such that $\Delta_{H_1}$ is an algebra homomorphism. \end{example}

One can also replace Hopf algebras in this paper by quasi-Hopf algebras for the corresponding non-strict versions. We leave this for future work. Beyond this we would expect by these means a full categorification where we let ${}_H\CM$ be replaced by any monoidal category $\CC$ with duals so that $\CZ(\CC)$ is braided, and  `categorify' by expressing all the data  above in categorical terms. We have set this up with Proposition~3.4 which gives the data in terms of the category and a braided Hopf algebra in the category $\CZ(\CC)$. Again, this is to be explored elsewhere.

\end{document}